\documentclass[twocolumn]{arxiv}

\usepackage{natbib}
\usepackage{amssymb}
\usepackage{amsmath}
\usepackage{graphicx}
\graphicspath{ {fig/} }
\usepackage{epstopdf}

\usepackage{float}

\floatstyle{ruled}
\newfloat{algorithm}{thp}{lop}
\floatname{algorithm}{Algorithm}
\newfloat{subroutine}{thp}{lop}
\floatname{subroutine}{Subroutine}

\newtheorem{assumption}{Assumption}
\newtheorem{system}{Dynamical system}

\begin{document}

\begin{frontmatter}
\title{Dynamic controllers for column synchronization of rotation matrices: a QR-factorization approach\thanksref{footnoteinfo}}

\thanks[footnoteinfo]{The authors gratefully acknowledge the financial support
from Fonds National de la Recherche Luxembourg (FNR8864515).}

\author[Johan]{Johan Thunberg}\ead{johan.thunberg@uni.lu}, 
\author[Johan]{Johan Markdahl}\ead{johan.markdahl@uni.lu}, 
\author[Johan]{Jorge Goncalves}\ead{jorge.goncalves@uni.lu}

\address[Johan]{Luxembourg Centre for Systems Biomedicine, University of Luxembourg,
    6, avenue du Swing, 
L-4367 Belvaux, Luxembourg}


\begin{abstract}
In the multi-agent systems setting, this paper addresses continuous-time distributed synchronization of columns of rotation matrices. More precisely, $k$ specific columns shall be synchronized and only the corresponding $k$ columns of the relative rotations between the agents are assumed to be available for the control design. When one specific column is considered, the problem is equivalent to synchronization on the $(d-1)$-dimensional unit sphere and when all the columns are considered, the problem is equivalent to synchronization on $\mathsf{SO}(d)$. We design dynamic control laws for these synchronization problems. The control laws are based on the introduction of auxiliary variables in combination with a QR-factorization approach. The benefit of this QR-factorization approach is that we can decouple the dynamics for the $k$ columns from the remaining $d-k$ ones. Under the control scheme, the closed loop system achieves almost global convergence to synchronization for quasi-strong interaction graph topologies. 
\end{abstract}

\end{frontmatter}

\section{Introduction}\label{sec:introduction}
This paper considers 
multi-agent systems continuously evolving on $\mathsf{SO}(d)$,
i.e., the set of $d \times d$ rotation matrices. The agents interact locally with each other and the neighborhood structure is determined by an interaction graph that is quasi-strongly connected. 
For such systems,
we address the following synchronization problem. How to design control laws in the body fixed coordinate frames of the agents such that $k$ specific columns of the rotation matrices asymptotically synchronize (converge to the set where they are the same and equal to the columns of a constant matrix) as time goes to infinity. The problem is, in general, a synchronization problem on a Stiefel manifold.
The control laws shall be designed by using the corresponding $k$ columns of the relative rotations between the agents (and not the other columns). Such control laws can be used to solve the synchronization problem on the unit sphere; consider for example the case where satellites in space only monitor one axis of each of its neighbors. But it can also be used in problems where various degrees of reduced attitudes are available, or the problem where complete rotations are available.  
To solve the problem we introduce auxiliary variables and use a QR-factorization approach. The benefit of this approach is that the dynamics of the $k$ columns considered can be decoupled from the dynamics of the remaining $d-k$ ones.

Two important special cases of the problem considered are synchronization of whole rotation matrices, i.e.,  synchronization on $\mathsf{SO}(d)$, and synchronization of one specific column vector, i.e., synchronization on the $(d-1)$-sphere. In these cases, for obvious reasons of applicability, the dimensions $d = 2$ and $d = 3$ have been mostly considered. The distributed synchronization problem on the unit sphere has been studied from various aspects~\cite{Sarlette2009,Olfati-Saber2006,Li2014,li2015collective}. Recently there have been some new developments~\cite{markdahl2016,markdahl2016towards}. 
\cite{markdahl2017almost,markdahl2016,markdahl2016global,pereira2016,lageman2016consensus}. 
In \cite{markdahl2017almost} the classical geodesic control law is studied for undirected graph topologies. Each agent moves in the tangent space in a weighted average of the directions to its neighbors. Almost global synchronization is de facto shown by a characterization of all the equilibria; the equlibria that are not in the synchronization set are shown to be unstable and the equilibra in the synchronization set are shown to be stable. The analysis can be seen to parallel the one in \cite{Tron2012} (also for undirected topologies) for the case of synchronization on $\mathsf{SO}(3)$, where intrinsic control laws are designed for almost global synchronization. For the case $d = 2$ an almost global synchronization approach has been presented for directed topologies and the $1$-sphere \cite{scardovi2007synchronization}. That approach is a special case of the one in \cite{AS-RS:09}.

The problem of synchronization on $\mathsf{SO}(3)$ has been extensively studied~\cite{Ren2010,Sarlette2010,Tron2013,Tron2014,thunbergaut,Deng2016}. Often the control algorithms are of gradient descent types and assume undirected topologies~\cite{distthunberg,sarlette2009automatica}. Local convergence results are often obtained~\cite{thunbergaut,johan02}. If a global reference frame is used, one can show almost global convergence~\cite{thunbergaut}---this is not allowed in the design of our control laws, only relative information is to be used. As mentioned above, \cite{Tron2012} provides a control algorithm for almost global convergence. The idea is to use so-called shaping functions where a gain constant can be chosen large enough to guarantee almost global consenus. The algorithm is defined in discrete time. 

By introducing auxiliary state variables based on the QR-factorization of matrices, this work provides a dynamic feedback control algorithm for synchronization of the $k$ first columns of the rotation matrices of the agents. The dynamics of the auxiliary variables follow a standard consensus protocol. The idea of using auxiliary or estimation variables with such dynamics is not new. Early works include \cite{scardovi2007synchronization} and \cite{AS-RS:09}, where the former addresses the $1$-sphere and the latter addresses manifolds whose elements have constant norms and satisfy a certain optimality condition. Such manifolds are $\mathsf{SO}(d)$ and the Grassmann manifold~$\text{Grass}(k,d)$. If the approach in \cite{AS-RS:09} is used to synchronize $k$ columns of the rotation matrix where $k < d -1$, then the entire relative rotations are used in the control design, which is not in general allowed in the problem considered here. In our proposed QR-factorization approach, only the corresponding $k$ columns of the relative rotations are used in the controllers. 
  Under the control scheme, the closed loop dynamics achieves almost global convergence to the synchronization set for quasi-strong interaction topologies.

\section{Preliminaries}\label{sec:preliminaries}
We start this section with some set-definitions.
We define the special orthogonal group 
$$\mathsf{SO}(d) = \{Q \in \mathbb{R}^{d \times d} : Q^TQ = I_d, \text{det}(Q) = 1\}$$
and set of skew symmetric matrices
$$\mathsf{so}(d) = \{\Omega \in \mathbb{R}^{d \times d} :  \Omega^T = - \Omega \}.$$
The $d$-dimensional unit sphere is 
$$\mathbb{S}^d = \{y \in \mathbb{R}^{d+1}: \|y\|_2 = 1\}.$$
The set of invertible matrices in $\mathbb{R}^{d \times d}$ is $$\mathsf{GL}(d) = \{Q \in \mathbb{R}^{d \times d} : \text{det}(Q) \neq 0\}.$$

We will make use of directed graphs, which have node set $\mathcal{V} = \{1,2,\ldots, n\}$
and edge sets $\mathcal{E} \subset \mathcal{V} \times \mathcal{V}$. Such a directed graph $\mathcal{G} = (\mathcal{V}, \mathcal{E})$ is quasi-strongly connected 
if it contains a rooted spanning tree or a center, i.e.,
there is one node to which there is a directed path 
from any other node in the graph. 
A directed path 
is a sequence of (not more than $n$) nodes such that any two consecutive nodes in the path comprises an edge in the graph. For $\mathcal{G} = (\mathcal{V}, \mathcal{E})$ we define 
$\mathcal{N}_i = \{j \in \mathcal{V}: (i,j) \in \mathcal{E}\}$ for all $i$. 

We will consider a multi-agent system with $n$ agents.
There are $n$ coordinate systems $\mathcal{F}_i$, each of which corresponding to a unique agent $i$ in the system. There is also a world (or global) coordinate system 
$\mathcal{F}_W$. At each time $t$, each coordinate system $\mathcal{F}_i$ is related to the global coordinate system $\mathcal{F}_W$ via a rotation $Q_i(t) \in \mathsf{SO}(d)$. This means that $Q_i(t)$ transforms vectors in $\mathcal{F}_i$ to vectors in $\mathcal{F}_W$. 

For all $i$, let $Q_i(t,k)$ be the ``tall matrix'' consisting of the first $k$ columns of 
$Q_i(t)$. Thus, $Q_i(t,d) = Q_i(t)$ and $Q_i(t,1)$ is the first column of $Q_i(t)$. All the columns of $Q_i(t,k)$ are obviously mutually orthogonal and each one an element of the $(d-1)$-sphere. Let $Q_{ij}(t) = Q_i^T(t)Q_j(t)$ and $Q_{ij}(t,k) = Q_i^T(t,d)Q_j(t,k)$ for all $i,j$. These matrices comprise the relative transformations between the coordinate frames $\mathcal{F}_j$ and $\mathcal{F}_i$ and the $k$ first columns thereof, respectively.

The matrix $R_i(t,d)$, or shorthand $R_i(t)$, is an element of $\mathbb{R}^{d \times d}$ for all $i$, $t$. The matrix $R_i(t,k) \in\mathbb{R}^{k \times k}$ is the upper left $k \times k$ block matrix of the matrix $R_i(t)$. These $R_i(t,k)$'s are communicated between the agents. For $R_i(t,k)$ invertible we define $R_{ij}(t,k)$ as $R_i(t,k)R_j^{-1}(t,k)$. Observe the difference in terms of the matrix inverse between $Q_{ij}(t,k)$ and  $R_{ij}(t,k)$, i.e.,  $Q_{ij}(t,k) = Q_{i}^{-1}(t,d)Q_j(t,k)$, whereas $R_{ij}(t,k) = R_{i}(t,k)R_j^{-1}(t,k)$.
Let $Q(t,k) = [Q_1^T(t,k),  Q_2^T(t,k), \ldots, Q_n^T(t,k)]^T \in \mathbb{R}^{nd \times k}$ and $R(t,k) = [R_1^T(t,k),  R_2^T(t,k), \ldots, R_n^T(t,k)]^T \in \mathbb{R}^{nk \times k}$.

The functions $\text{low}(\cdot)$ and $\text{up}(\cdot)$ are defined for matrices in $\mathbb{R}^{m_1 \times m_2}$ for all $m_1 \geq m_2$. The function $\text{low}(\cdot)$ returns a matrix of the same dimension as the input, a matrix in $\mathbb{R}^{m_1 \times m_2}$ that is, where 
each $(i,j)$-element of the matrix is equal to that of the input matrix if $i > j$ and equal to $0$ if $i \leq j$. The function $\text{up}(\cdot)$ returns a matrix in $\mathbb{R}^{m_2 \times m_2}$; each $(i,j)$-element of the matrix is equal to that of the input matrix if $i                                  \leq j$ and equal to $0$ if $i > j$. 

We continue by introducing two assumptions that will be used in the problem formulation in the next section.

\begin{assumption}[Connectivity]\label{ass:1}
It holds that $\mathcal{G}$ is 
quasi-strongly connected.
\end{assumption}

\begin{assumption}[Dynamics]\label{ass:2}
The time evolution of the state of each agent $i$ is given by
\begin{align}\label{eq:v3:501}
\frac{d}{dt}{Q}_i(t,d)  = Q_i(t,d)U_{i}(t,d), 
\end{align}
where $U_i(t,d) \in \mathsf{so}(d)$ and $Q_i(0,d) \in \mathsf{SO}(d)$. In particular it holds that 
\begin{align}\label{eq:v3:1}
\frac{d}{dt}{Q}_i(t,k)  = Q_i(t,d)U_{i}(t,k), \forall  k \in \{ 1,2, \ldots, d\},
\end{align}
where $U_i(t,k) =  U_i(t,d)[I_k, 0]^T$. 
\end{assumption}
The $U_i(t,d)$'s are the controllers we are to design. 
An important thing to note in \eqref{eq:v3:501} is that $U_i(t,d)$, or rather the columns thereof, are defined in the $\mathcal{F}_i$-frames. If those would have been defined in the world frame $\mathcal{F}_{W}$, the agents would have needed to know their own rotations to that frame, i.e., the $Q_i$-matrices. Those matrices are not assumed to be available for the agents.

We let $(\mathsf{SO}(d))^n$ be the following subset of $\mathbb{R}^{nd \times d}$,
\begin{align*}
\{Z : Z = [Z_1^T, Z_2^T, \ldots, Z_n^T]^T, Z_i \in \mathsf{SO}(d)\:  \forall i \}.
\end{align*}
We let $(\mathsf{GL}(d))^n$ be the following subset of $\mathbb{R}^{nd \times d}$,
\begin{align*}
\{Z : Z = [Z_1^T, Z_2^T, \ldots, Z_n^T]^T, Z_i \in \mathsf{GL}(d) \:  \forall i \}. 
\end{align*}

\section{Problem formulation}
The goal is to design $U_i(t,d)$ (and in particular $U_i(t,k)$) as a dynamic feedback control law such that the $Q_i(t,k)$-matrices asymptotically aggregate or converge to the synchronization set. In this problem, formally presented below, a key assumption in the control design is that the information available for the agents is both relative and local. This means that no knowledge of a global coordinate system is assumed and that only relative rotations and vectors between local neighboring agents are used.

\begin{prob}\label{prob:1}
Let Assumption~\ref{ass:1} and Assumption~\ref{ass:2} hold. Let $d \geq 2$ and $k \leq d -1$.
Design the $U_i(t,d)$-controllers as continuous functions of $t$ and the elements in the collection $\{Q_{ij}(t,k)\}_{j \in \mathcal{N}_i}$ such that the following is fulfilled: there is a unique continuous solution for the $Q_i$'s and
there is $\bar{Q} \in \mathsf{SO}(d)$ such that
$$\lim_{t \rightarrow \infty}\|Q_i(t,k) - \bar{Q}[I_k, 0]^T\| = 0~\forall i.$$
\end{prob}
In Problem~\ref{prob:1} we made the assumption that $k \leq d -1$. 
This  restriction of the $k$'s to those smaller than or equal to $d-1$ can be made without loss of generality. If all the $Q_i(t,d-1)$'s are equal, so are all
the $Q_i(t,d)$'s.
This means that synchronization of the $Q_i(t,d-1)$'s is equivalent to synchronization of the $Q_i(t,d)$'s. 

For $k = 1$, our problem was studied in e.g., \cite{markdahl2016,markdahl2016global,pereira2016} and for $k = 1$, $d = 2$ in \cite{scardovi2007synchronization}. For the case $k = d$, the problem has been studied in e.g.,  \cite{thunbergaut,Deng2016,Ren2010,Sarlette2010,Tron2013,Tron2014}.
In general, for $k < d-1$, Problem~\ref{prob:1} is related to  synchronization on the Grassmann manifold~$\text{Grass}(k,d)$, which was studied in~\cite{AS-RS:09}. It is de facto synchronization on a compact Stiefel manifold.

\section{The proposed control algorithm}\label{sec:4}
In this section we propose a control algorithm as a candidate solution to Problem~\ref{prob:1}. 

Before we  introduce the control algorithm, Algorithm 1 below, we define the set
\begin{align*}
\mathcal{D}_{QR}(k) & = (\mathsf{SO}(d)[I_k,0]^T)^n \times (\mathcal{R}^+(k))^n.
\end{align*}
The set $\mathcal{R}^+(k)$ comprise the upper triangular matrices in $\mathbb{R}^{k \times k}$, whose diagonal elements are positive. 
The algorithm below is restricted to those $[Q(t,k)$, $R(t,k)]$'s contained in $\mathcal{D}_{QR}(k)$. 

\begin{algorithm}[h!]
\caption{Distributed control algorithm for synchronization of the $Q_i(t,k)$'s}\label{alg1}
\vspace{2mm}
\textbf{Initialization at time $0$:} for all $i$, choose $d \geq 2$ and $k \leq d -1$, let $R_i(0,k)$ be upper triangular with positive elements on the diagonal, and let $Q_i(0,d) \in {SO}(d)$.\\  \\
\textbf{Inputs to agent $i$ at time $t \geq 0$:}  $(Q_{ij}(t,k))_{j \in \mathcal{N}_i}$ and  $(R_{j}(t,k))_{j \in \mathcal{N}_i}$. \\ \\
\textbf{Controllers at time $t \geq 0$ when $R_i(t,k) \in \mathsf{GL}(k)$:} \\ \\let, for all $i$, 
\begin{align}
\label{eq:main:1}
V_i(t,k) & = \sum_{j \in \mathcal{N}_i}a_{ij}(Q_{ij}(t,k)R_{ji}(t,k) - [I_k, 0]^T), \\
\label{eq:main:4}
U_i(t,d) & = [\text{low}(V_i(t,k)), 0] - [\text{low}(V_i(t,k)), 0]^T \\
\label{eq:main:3}
\dot{R}_i(t,k) & = \text{up}((V_i(t,k) - U_i(t,k))R_i(t,k)), \\
\nonumber
& \quad \text{ where } \\
\label{eq:main:2}
U_i(t,k) & = \text{low}(V_i(t,k)) - [I_k, 0]^T(\text{low}(V_i(t,k))^T[I_k, 0]^T), \\
\nonumber
& \quad \text{ comprises the first k columns of } U_i(t,d). 
\end{align}
\end{algorithm}

We will show in the next section that the restriction of the $[Q(t,k)$, $R(t,k)]$'s to those contained in $\mathcal{D}_{QR}(k)$ does not comprise a limitation from a practical point of view since the set $\mathcal{D}_{QR}(k)$ is invariant under the proposed control scheme for all but a set of measure zero of the initial points.

It is important to note in Algorithm 1 (see next page) that $U_i(t,d)$ (and not $U_i(t,k)$) is the controller that is used by agent $i$. The latter is however the part of the controller that affects $Q_{i}(t,k)$ according to \eqref{eq:v3:1}. We remind the reader that $U_i(t,k)$ comprises the first $k$ columns of $U_i(t,d)$, i.e., it is a restriction of $U_i(t,d)$.  At the initialization step, the matrix $R_i(0,k)$ is chosen (or constructed) by each agent $i$, whereas the matrix $Q_i(0,k)$ is not. The latter is not known by the agent under our ``relative information only''-assumption in the control design.

The expressions \eqref{eq:main:1}-\eqref{eq:main:2} in Algorithm 1 seem, at a first glance, a bit non-intuitive. As it turns out, the closed-loop system under Algorithm 1 is, almost everywhere, equivalent to the QR-decompositions of the $Z_i$-matrix variables in Dynamical System~\ref{sys:1} below, which is a standard linear continuous time consensus protocol. In the next section we will, through a series of technical results, prove the equivalence between Algorithm 1 and Dynamical System~\ref{sys:1} (almost everywhere).  

\begin{system}\label{sys:1}
Let Assumption~\ref{ass:1} hold and let $k \leq d-1$. 
Let $Z(t,k) = [Z_1(t,k), Z_2(t,k),
\ldots, Z_n(t,k)]$ for all $t \geq 0$, where the $Z_i(t,k)$'s (including the initial points at time $0$) are elements in $\mathbb{R}^{d \times k}$ for all $t \geq 0$. The time evolution of $Z(t,k)$ is governed by   
\begin{align}\label{eq:13}
\dot{Z}_i(t,k) & = \sum_{j \in \mathcal{N}_i}a_{ij}(Z_j(t,k) - Z_i(t,k))~ \forall i, \text{ where } \\
\nonumber
& \quad ~ a_{ij} > 0 \: \forall (i,j).
\end{align}
\end{system}

\section{Convergence results for Algorithm 1}
Suppose Assumption~\ref{ass:1} and Assumption~\ref{ass:2} hold 
and the signals $Q_i(t,k)$ and $R_i(t,k)$
are produced by Algorithm~1. We define 
\begin{align*}
& \mathcal{D}_{QR, \text{full}}(k) \\ 
=~& \{[Q(0,k), R(0,k)]:  \\
& 1) \hspace{1mm}~[Q(t,k), R(t,k)] \in \mathcal{D}_{QR}(k) ~\forall t \geq 0, \\ 
& 2) \hspace{1mm}~\exists~ \bar{Q} \in \mathsf{SO}(d) \text{ and } \bar{R} \in \mathsf{GL}(k), \text{ s.t. } \\
& \lim_{t \rightarrow \infty}[Q_i(t,k), R_i(t,k)] = [\bar{Q}[I_k, 0]^T, \bar{R}]~\forall i \}.
\end{align*}
A verbal interpretation of the set $\mathcal{D}_{QR, \text{full}}(k)$ is the following one. It is the set of initial points in $\mathcal{D}_{QR}(k)$ such that 1) all the $R_i(t,k)$'s are invertible at all times $t$ larger than $0$, and 2) the $Q_i(t,k)$'s and the $R_i(t,k)$'s converge to matrices $\bar Q [I_k, 0]^T$ respective $\bar{R}$, where $\bar Q$ is in $SO(d)$ and  $\bar{R}$ is invertible. 

Now, it would be good if we could prove that the set $\mathcal{D}_{QR, \text{full}}(k)$ contains most of $\mathcal{D}_{QR}(k)$.  
The following proposition provides such a result; it is the main result of the paper. It is claiming almost global convergence to the synchronization set for Algorithm 1. For all but a set of measure zero of $Q_i(0,k)$'s and $R_i(0,k)$'s, the matrices converge to the synchronization set. 
The rest of this section is dedicated to the proof of the claim in the proposition.

\begin{prop}\label{prop:1}
Suppose Assumption~\ref{ass:1} and Assumption~\ref{ass:2} hold and the signals $Q_i(t,k)$ and $R_i(t,k)$
are produced by Algorithm 1. The set $\mathcal{D}_{QR, \text{full}}(k)$ is open and 
the set $\mathcal{D}_{QR}(k) - \mathcal{D}_{QR, \text{full}}(k)$ has measure zero and is nowhere dense in $\mathcal{D}_{QR}(k)$. 
\end{prop}

By inspection, we can verify that if all the $Q_i(t,k)$'s are equal and all the $R_i(t,k)$'s are equal, then all the expressions in \eqref{eq:main:1}-\eqref{eq:main:2} are equal to zero, which means that the system is at equilibrium. Now, what we want to establish is the almost global convergence to such an equilibrium. However, the structure of the system seems at first hand complicated, which might make the convergence analysis cumbersome. 

Now, instead of studying the dynamics of the closed loop system under \eqref{eq:main:1}-\eqref{eq:main:2}, 
the main idea of the proof of Proposition~\ref{prop:1} is to show that (as already mentioned in Section~\ref{sec:4}) when $[Q(0,k),R(0,k)]$ is contained in $\mathcal{D}_{QR, \text{full}}(k)$, there is a change of coordinates so that after this change of coordinates the dynamics of the system is described by the simple consensus protocol in Dynamical system 1 for the $Z_i$-matrices in $\mathbb{R}^{d \times k}$, see \eqref{eq:13} below. The idea is to, instead of studying the convergence of a seemingly complicated dynamical system, use well-known results for the simple consensus protocol. This type of approach is indeed not new, see for example \cite{scardovi2007synchronization} and \cite{AS-RS:09}. 

A complication here is that in order to establish the convergence result in Proposition~\ref{prop:1} we need to guarantee that in general the matrices evolving under the consensus protocol have full rank for all time points along the trajectories (as well as in the limit). This result is provided by Lemma~\ref{lem:2} below. After the lemma (and its proof) has been provided, we give the proof of Proposition~\ref{prop:1}.

We begin by recalling the following known result.

\begin{prop}[\cite{dieci1999smooth}]\label{prop:2}
Any full column rank time-varying $\mathcal{C}^k$-matrix has a $\mathcal{C}^k$ QR-decomposition. 
\end{prop}

We define 
\begin{align*}
\mathcal{D}_{Z}(k) ={} & \{X = [X_1^T, X_2^T, \ldots, X_n^T]^T \in \mathbb{R}^{nd \times k}: \\
& \hspace*{1mm}~ X_i \in \mathbb{R}^{d \times k} \text{ has full column rank } \forall i, \}
\end{align*} 
and define the set $\mathcal{D}_{Z, \text{full}}(k)$ as those $Z(0,k)$'s in $\mathcal{D}_{Z}(k)$ for which it holds that  $Z(t,k) \in \mathcal{D}_{Z}(k)$ 
for all $t$, when $Z(t,k)$ is generated by Dynamical System~\ref{sys:1} and there exists $\bar Z$  that has full column rank such that $\bar Z = \text{lim}_{t\rightarrow \infty}Z(t,k)$. 

\begin{lem}\label{lem:2}
For Dynamical System~\ref{sys:1} the following statements hold
\begin{enumerate}
\item there is a unique analytic solution for \eqref{eq:13}; \\
\item let $H(t,k)$ denote the convex hull of the $Z_i(t,k)$'s. The set $(H(t,k))^n$ is forward invariant for any $t$ and there is $\bar Z \in \mathbb{R}^{d \times k}$ (as a function of the initial state) such that $\|Z_i(t,k) - \bar{Z}\|_{\text{F}}$
goes to zero as $t$ goes to infinity;\\
\item for any time interval $[0, t_1)$ during which the $Z_i(t,k)$'s have full column rank, if (instead of being produced by Algorithm 1) the $Q_i(t,k)$'s and the $R_i(t,k)$'s correspond to QR-decompositions of the $Z_i(t,k)$'s (i.e., $Z_i(t,k) = Q_i(t,k)R_i(t,k)$ for all $i$),  then the $Q_i(t,k)$'s and the $R_i(t,k)$'s can be chosen as smooth functions (of $t$), where $Q_i(t,d) \in \mathsf{SO}(d)$, and $R_i(t,k) \in \mathsf{GL}(k)$ for all $i$;  \\
\item for all but a set of measure zero of the initial points $Z(0,k) \in \mathcal{D}_{Z}(k)$, the matrix $\bar{Z}$ in (2) has full column rank; \\
\item $\mathcal{D}_{Z}(k)-\mathcal{D}_{Z, \text{full}}(k)$ has measure zero and is nowhere dense. The set $\mathcal{D}_{Z, \text{full}}(k)$ is open in 
$\mathbb{R}^{nd \times d}$.
\end{enumerate}
\end{lem}

\quad \emph{Proof}: 
\emph{(1)} The dynamics for $Z$ can be written as 
\begin{align}
\label{eq:olle:1}
\dot{Z} = -(L \otimes I_d)Z,
\end{align}
where $L$ is the weighted graph Laplacian matrix defined by: $[L]_{ij} = 0$ if $i \neq j$ and $(i,j) \not\in \mathcal{E}$;  $[L]_{ij} = -a_{ij}$ if $i \neq j$ and $(i,j) \not\in \mathcal{E}$; $[L]_{ii} = \sum_{j \in \mathcal{N}_i}a_{ij}$. The matrix $L$ is implicitly parameterized by the $a_{ij}$'s. The unique analytic solution to \eqref{eq:olle:1} is 
\begin{equation}
Z(t,k) = A(t)Z(0,k),
\end{equation}
where $A(t) = \text{exp}(-(L \otimes I_d)t)$.  

\emph{(2)} \eqref{eq:olle:1} is convergent due to the properties of the graph Laplacian matrix $L$. A proof of statement \emph{(2)} is obtained by direct application of the results in Z. Lin et al (2007) or \cite{shi2009global}.

\emph{(3)} Direct application of \emph{(1)} and Proposition~\ref{prop:2}. We can without loss of generality assume that the $R_i(t,k)$'s have positive elements on the diagonal. Suppose it was not the case, then we multiply $Q_i(t,k)$ and $R_i(t,k)$ with a diagonal matrix from the right and the left, respectively. This diagonal matrix is constant and has $1$'s at $ii$-entries where $[R_i(t,k)]_{ii}$ is positive and $-1$'s at the $ii$-entries where $[R_i(t,k)]_{ii}$ is negative. Since $k \leq d-1$, we can choose this diagonal matrix in such a way that the resulting  orthogonal matrix after the multiplication (from the right) with the diagonal matrix is an element of $\mathsf{SO}(d)$. 

{\emph{(4)}} 
\eqref{eq:olle:1} is translation invariant. Translation invariance means that if we disturb the initial $Z_i(0,k)$'s by adding a matrix $\Xi$ to all the $Z_i(0,k)$'s, the matrix $\Xi$ is canceled out in the dynamics \eqref{eq:olle:1}, i.e., it is invariant to this common translation of the initial states. As a consequence, the state trajectories for the disturbed initial conditions are equal to those without the disturbance up to the added $\Xi$. This also holds in the limit so the equivalent matrix to $\bar Z$ for such disturbed initial conditions is $\bar Z + \Xi$.

It holds that $\mathbb{R}^{nd \times k} - \mathcal{D}_{Z}(k)$ has measure zero. The rest of this proof is about showing that the subset of $\mathcal{D}_{Z}(k)$ under consideration comprises all but a measure zero of $\mathbb{R}^{nd \times k}$. 
Now we define 
\begin{align*}
\mathcal{A}  ={}& \{X = [X_1^T, X_2^T, \ldots, X_n^T]^T \in \mathbb{R}^{nd \times k}: \\
& \hspace*{1mm}~ X_i = X_j \in \mathbb{R}^{d \times k} ~\forall i,j \}.
\end{align*}
In the above, $\mathcal{A}$ and $\mathcal{A}^{\perp}$ are linear subspaces of $\mathbb{R}^{nd \times k}$. The set $\mathcal{A}$ is the synchronization set and the set $\mathcal{A}^{\perp}$ is its orthogonal complement defined via the standard trace inner product. Each $Z(0,k)$ can be written as a sum of two matrices $X = [X_1^T, X_1^T, \ldots, X_1^T]^T$ and $Y = [Y_1^T, Y_2^T, \ldots, Y_n^T]^T$, where $X$ is an element of $\mathcal{A}$ and $Y$ an element of $\mathcal{A}^{\perp}$. Furthermore it holds that $\bar{Z}(X + Y) = \bar{Z}(Y) + X_1,$ where $\bar{Z}$ is given as a function of the initial condition. Now, for any fixed $Y \in \mathcal{A}^{\perp}$ it holds that $\bar{Z}(X + Y) = \bar{Z}(Y) + X_1$ has full column rank for all but a set of measure zero of the $X_1$'s in $\mathbb{R}^{d \times k}$. 

{\emph{(5)}} We begin by proving the zero measure of $\mathcal{D}_{Z}(k)-\mathcal{D}_{Z, \text{full}}(k)$.
The columns of $A(t)$ are linearly independent for all $t$ and analytic in $t$. 
Thus the matrix $A^T(t)$ has a smooth QR-factorization $\bar{Q}(t)\bar{R}(t)$, see Proposition~\ref{prop:2}.

We will without loss of generality only consider the case
$i = 1$ and prove that $Z_1(t,k)$ will be of full rank for all $t \geq 0$ for all but a set of measure zero of the $Z(0,k)$'s in $\mathcal{D}_{Z}(k)$. We can make this ``wlog-assumption'' since the procedure of the proof is equivalent for the other choices of $i$ after a permutation of the $Z_i$-matrices.

Let $z^1(t)$ be the first column of $Z(t,k)$, $z^2(t)$ the second column and so on. Let $Z_1(t,k) = B(t)Z(0,k)$, where 
\begin{align*}
B(t) & = [I_d, 0, \ldots, 0]A(t), 
\end{align*} 
Now we define
\begin{align*}
C(t,0) & = ([0,I_{(n-1)d}]\bar Q^T(t))^T, \\
C(t,1) & = [([0,I_{(n-1)d}]\bar Q^T(t))^T, z^1(0)], \\
& \: \: \:\vdots \\
C(t,d) & = [([0,I_{(n-1)d}]\bar Q^T(t))^T, z^1(0),\ldots, z^d(0)].
\end{align*} 
The $C$'s (besides $C(t,0)$) are implicitly parameterized by the $Z(0,k)$'s.
It holds that $\text{im}(C(t,0)) = \text{ker}(B(t))$ for all $t$. The columns of $C(t,0)$ comprise a smooth orthogonal basis for $\text{ker}(B(t))$. The  $C(t,k)$'s for $k = 1,2, \ldots, d$ are built in the proposed way to take into account that, for a singular $Z_i$, one of
the columns of a $C(t,k)$-matrix is a linear combination of previous columns in the matrix. Formally, these matrices are used in the following condition: if there is a $t_1$ such that the rank of $Z_1(t_1,k)$ is smaller than $k$, there must be a $\bar{k} \in \{1,2,\ldots,k\}$ with a vector $\bar a \in \mathbb{R}^{(n-1)d + \bar{k}-1}$ such that $C(t_1,\bar{k}-1)\bar a = z^{\bar{k}}(0)$. We refer to this as \textbf{Condition} \textbf{a)}.

Now, the reminder of the proof amounts to showing that 
for $\bar k \in \{1,2, \ldots, k\}$ only for a set of measure zero of $z^k(0)$'s there is $[t,\bar a^T]^T \in \mathbb{R}^{(n-1)d + \bar{k}}$ such that 
Condition a) is fulfilled. Once this is showed we can conclude that $Z(t,{k})$ has full rank for all $t$ for all but a set of measure zero of the initial conditions (i.e.,  $Z(0,{k})$). If the latter would not have been true, there would have been a $\bar{k}$ for which Condition a) is satisfied for a positive measure of $z^{\bar{k}}(0)$'s, but that would have been a contradiction. 

For all $\bar k$ we define 
$$g_{\bar{k}}(t,a) = C(t,\bar{k}-1)a,$$
which is a smooth mapping from $\mathbb{R}^{(n-1)d + \bar{k}}$ to 
$\mathbb{R}^{nd}$ (we allow for negative $t$'s here).
It is a nonlinear (and smooth) function of $t$ and a linear function of $a$. The elements of $a$ comprise weights in the sum of columns of $C(t,\bar{k}-1)$, which is the returned vector of $g_{\bar{k}}(t,a)$. We want to show that we cannot in general choose $t$ and $a$ such that this column sum is equal to $z^{\bar{k}}(0)$, thereby showing that Condition a) is not fulfilled in general. Now, the rank of the Jacobian matrix of $g_{\bar{k}}$ is at most $(n-1)d + \bar k < nd$ for all $[t,a^T]^T \in \mathbb{R}^{(n-1)d + \bar k}$. Hence all $[t,a^T]^T$ are critical points and, due to Sard's theorem, $\text{im}(g_{\bar{k}})$ has measure zero in $\mathbb{R}^{nd}$. This means that Condition a) is only fulfilled for a set of measure zero of $z^{\bar{k}}(0)$'s. 

We know that for all but a set of measure zero of the initial points $Z(0,k) \in \mathcal{D}_{Z}(k)$, the matrix $\bar{Z}$---the matrix that the states converge to---has full column rank. It holds that $\mathbb{R}^{nd \times k} - \mathcal{D}_{Z}(k)$ has measure zero. We know that Condition a) is only fulfilled for a set of measure zero of $z^{\bar{k}}(0)$'s, see above.
Thus we can conclude, $\mathcal{D}_{Z}(k)- \mathcal{D}_{Z, \text{full}}(k)$ has measure zero.

Now, we prove that $\mathcal{D}_{Z}(k)- \mathcal{D}_{Z, \text{full}}(k)$ is nowhere dense in $\mathcal{D}_{Z}(k)$. We prove this by showing that for any neighborhood $U$ in $\mathcal{D}_{Z}(k)$, the set $(\mathcal{D}_{Z}(k)- \mathcal{D}_{Z, \text{full}}(k)) \cap U$ is not dense in $U$. Note that $\mathcal{D}_{Z}(k)$ is open in $\mathbb{R}^{nd \times k}$.

Consider an arbitrary neighborhood $U$ in $\mathcal{D}_{Z}(k)$. There must be a the point $Z_0 \in U \cap \mathcal{D}_{Z, \text{full}}(k)$, since the set $\mathcal{D}_{Z}(k)- \mathcal{D}_{Z, \text{full}}(k)$ has measure zero. Now, since  $Z_0 \in \mathcal{D}_{Z, \text{full}}(k)$, the following holds (each statement in the following list continues on the previous ones):
there is $\bar{Z} \in \mathcal{D}_{Z}(k)$ such that for $Z(t,k)$ initialized at the point $Z_0$ at time $0$, the $Z_i(t,k)$'s converge to $\bar{Z}$;  
there is a closed ball $B_{\bar Z,\epsilon} \subset \mathcal{D}_Z(k)$ with radius $\epsilon > 0$ centered at $\bar{Z}$; there is a finite time $t_f > 0$ after which $Z(t,k) \in B_{\bar Z,\epsilon/2}$, where $B_{\bar Z,\epsilon/2}$ is the closed ball around $\bar Z$ with radius $\epsilon/2$. The two balls above are defined with respect to the Frobenius norm.

Now, define $f(Z) = \inf_{Y \in \mathbb{R}^{nd \times k} - \mathcal{D}_{Z}(k)}\|Y - Z\|$, which is a continuous function on $\mathbb{R}^{nd \times k}$.  This means that the composite function $(f \circ Z)(t,k,Z_0)$ is continuous in $t$ on $[0,t_f]$ (we treat $k$ and the initial condition $Z_0$ as parameters). Furthermore, the function is strictly positive on $[0,t_f]$ and attains a minimum $\epsilon_2 >0$ there. Now choose $\epsilon_3$ such that $0 < \epsilon_3 < \min\{\epsilon/2, \epsilon_2\}$. Due to the continuous dependency theorem on initial conditions, there is $\epsilon_4 > 0$ such that if $Z_{0, \text{new}} \in B_{Z_0, \epsilon_4}$, then $\|(f \circ Z)(t,k,Z_0) - (f \circ Z)(t,k,Z_{0,\text{new}})\| \leq \epsilon_3$ for all $t\in [0,t_f]$. For any point $Z_{0,\text{new}} \in B_{Z_0, \epsilon_4}$ it holds that if $Z(t,k)$ is initialized at $Z_{0,\text{new}}$, then the $Z_i(t)$'s are contained $B_{\bar Z,\epsilon}$ after time $t_f$. 
The ball $B_{\bar Z,\epsilon}$ is forward invariant under \eqref{eq:13}
since the set $(H(t,k))^n$ is contained in $B_{\bar Z,\epsilon}$ and is forward invariant, see Lemma~\ref{lem:2} \textit{(2)}. Thus, for all times $t \geq t_f$, all the $Z_i(t,k)$'s have full rank (as well as for the limit matrix) when the initial condition for $Z(t,k)$ was chosen to $Z_{0,\text{new}}$.

What we have shown is that there is a ball $B_{Z_0,\epsilon_4}$ around $Z_0$ for which no point is contained in $(\mathcal{D}_{Z}(k)- \mathcal{D}_{Z, \text{full}}(k))$, i.e., it is an interior point of $\mathcal{D}_{Z, \text{full}}(k)$ as well as $\mathbb{R}^{nd \times k}$. This means that  $(\mathcal{D}_{Z}(k)- \mathcal{D}_{Z, \text{full}}(k))$ is not dense in $U$. But $U$ was arbitrarily chosen in $\mathcal{D}_{Z}(k)$. Thus the proof is complete. We also obtain the result that $\mathcal{D}_{Z, \text{full}}(k)$ is open in $\mathbb{R}^{nd \times k}$. 
\hfill $\blacksquare$

\begin{rem}
In the proof of \emph{(5)} when we prove that $\mathcal{D}_{Z}(k)- \mathcal{D}_{Z, \text{full}}(k)$ has measure zero, we did not use the property that the matrix in the right-hand side of \eqref{eq:olle:1} is $-(L \otimes I_d)$. That partial result also holds for any other matrix in $\mathbb{R}^{nd \times nd}$. Furthermore, a straightforward generalization is  that only for a measure zero set of the initial matrices in $\mathbb{R}^{nd \times k}$, any matrix block of $Z$ with more rows than columns loses rank at a finite time.
\end{rem}

\quad \emph{Proof of Proposition~\ref{prop:1}}: 
We begin by returning to the system Dynamical system 1. We suppose Assumption~\ref{ass:1} and Assumption~\ref{ass:2} hold and that $Z(t,k)$ is the solution of \eqref{eq:13}. 

There is a diffemorphism $h$ between $\mathcal{D}_{Z}(k)$ and $\mathcal{D}_{QR}(k)$. For every $X = [X_1^T, X_2^T, \ldots, X_n^T]^T$ in $\mathcal{D}_{Z}(k)$, $h$ provides the unique $[\tilde{Q}, \tilde{R}] = [[\tilde{Q}_1^T, \ldots, \tilde{Q}^T_n]^T,[\tilde{R}_1^T, \ldots, \tilde{R}_n^T]^T]$ in $\mathcal{D}_{QR}(k)$ such that $X_i = \tilde{Q}_i\tilde{R}_i$ for all $i$. To see that it is a diffeomorphism we note the following. Let $h(X_i) = [\tilde Q_i, \tilde R_i]$. The $\tilde R_i$-matrix is the upper triangular matrix obtained by Cholesky's factorization of $X_i^TX_i$. Cholesky's factorization is analytic on the set of positive definite matrices~\cite{lee2014matrix}. The $\tilde Q_i$-matrix is given by $\tilde Q_i = X_i\tilde R_i^{-1}$, where $\tilde R_i^{-1}$ is analytic on the set of invertible matrices. As for $h^{-1}$, it holds that $h^{-1}([\tilde Q_i, \tilde R_i]) = X_i=\tilde Q_i\tilde R_i$, where the right-hand side in the last equation is analytic in the two matrices. 
 
\begin{figure}[!th]
\centering
\includegraphics[scale=0.1]{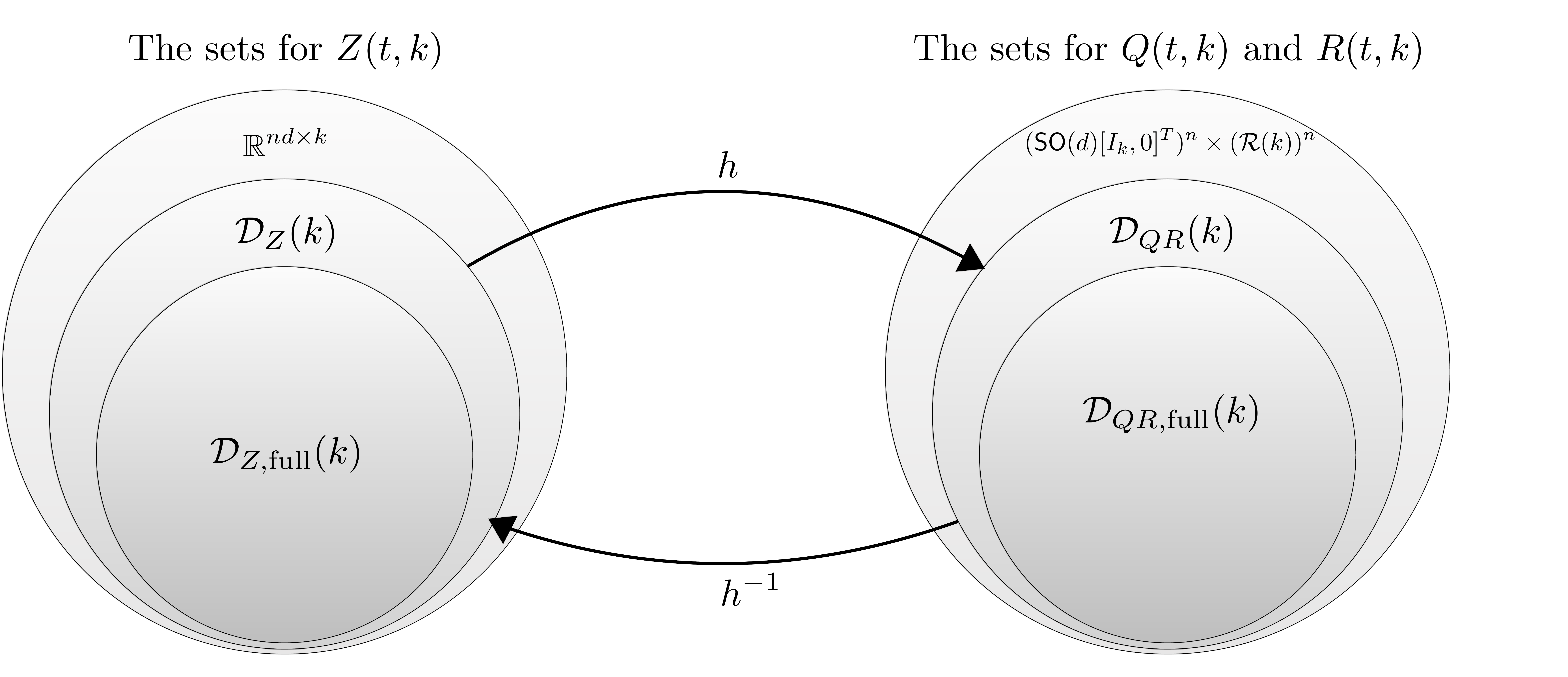}
\caption{Relations between the matrix sets.}
 \label{fig:1001}
\end{figure}

Figure~\ref{fig:1001} illustrates relations between the matrix sets. 

It holds that $\mathbb{R}^{nd \times k} - \mathcal{D}_{Z}(k)$ has measure zero. 
It also holds that set $\mathcal{D}_{Z}(k) - \mathcal{D}_{Z, \text{full}}(k)$ has measure zero, see (4) and (5) in Lemma~\ref{lem:2}.
The proof amounts to showing that $\mathcal{D}_{QR, \text{full}}(k) = h(\mathcal{D}_{Z, \text{full}}(k))$. 

We know that $\mathcal{D}_{Z}(k) - \mathcal{D}_{Z, \text{full}}(k)$ has measure zero and that a measure zero set is mapped to a measure zero set under a diffeomorphism. Hence, if we show that  $\mathcal{D}_{QR, \text{full}}(k) = h(\mathcal{D}_{Z, \text{full}}(k))$, then we know that $h(\mathcal{D}_{Z}(k) - \mathcal{D}_{Z, \text{full}}(k)) = \mathcal{D}_{QR}(t,k) - \mathcal{D}_{QR, \text{full}}(t,k)$ and that the set on the right-hand side has measure zero. Furthermore, open sets are mapped to open sets and $\mathcal{D}_{Z, \text{full}}(k)$ is open. Thus $\mathcal{D}_{QR, \text{full}}(k)$ is open if $\mathcal{D}_{QR, \text{full}}(k) = h(\mathcal{D}_{Z, \text{full}}(k))$. This latter fact in combination with the fact that $\mathcal{D}_{QR}(t,k) - \mathcal{D}_{QR, \text{full}}(t,k)$  has measure zero can be used to conclude that $\mathcal{D}_{QR}(t,k) - \mathcal{D}_{QR, \text{full}}(t,k)$ is nowhere dense.

We first show that \textbf{1)} $h(\mathcal{D}_{Z, \text{full}}(k)) \subset \mathcal{D}_{QR, \text{full}}(k)$ and then show that \textbf{2)} $h^{-1}(\mathcal{D}_{QR, \text{full}}(k)) \subset \mathcal{D}_{Z, \text{full}}(k)$.

\textbf{1)}
We assume that $Z(0,t)$ is in $\mathcal{D}_{Z, \text{full}}(k)$ and define
\begin{align*}
Z(t,k) & = [Z_1^T(t,k),\ldots, Z_n^T(t,k)]^T \\
& = [(\tilde Q_1(t,k) \tilde R_1(t,k))^T, \ldots, (\tilde Q_n(t,k)\tilde R_n(t,k))^T]^T,
\end{align*}
where $[\tilde{Q}_i(t,k), \tilde{R}_i(t,k)]$ is the unique QR-decomposition of $Z_i(t,k)$ with positive diagonal elements for $\tilde{R}_i(t,k)$. These QR-decompositions are smooth according to Proposition~\ref{prop:2}. 

Now we prove that the dynamics for the $\tilde{Q}_i(t,k)$'s and $\tilde{R}_i(t,k)$'s have the same structure as \eqref{eq:main:1}-\eqref{eq:main:3}.
For all $i$, it holds that 
\begin{align*}
\dot{Z}_i(t,k) & =  \dot{\tilde{Q}}_i(t,k)\tilde R_i(t,k) + \tilde{Q}_i(t,k)\dot{\tilde{R}}_i(t,k) \\
& = \tilde{Q}_i(t,d)\tilde{U}(t,k)\tilde{R}_i(t,k) + \tilde{Q}_i(t,k)\dot{\tilde{R}}_i(t,k) 
\end{align*}
and it also holds that 
\begin{align*}
\dot{Z}_i(t,k) & = \sum_{j \in \mathcal{N}_i}a_{ij}(Z_j(t,k) - Z_i(t,k)) \\
& = \sum_{j \in \mathcal{N}_i}a_{ij}(\tilde{Q}_j(t,k)\tilde R_j(t,k) - \tilde{Q}_i(t,k)\tilde R_i(t,k)).
\end{align*}
Thus, 
\begin{align*}
& \quad ~ \tilde{Q}_i(t,d)\tilde{U}(t,k)\tilde{R}_i(t,k) + \tilde{Q}_i(t,k)\dot{\tilde{R}}_i(t,k) \\
& = \sum_{j \in \mathcal{N}_i}a_{ij}(\tilde{Q}_j(t,k)\tilde R_j(t,k) - \tilde{Q}_i(t,k)\tilde R_i(t,k)).
\end{align*}
We can express the last equation above as 
\begin{align}
\label{eq:niss:1}
& \quad \: \: \tilde{Q}_i(t,d)\tilde{U}(t,k) + \tilde{Q}_i(t,k)\dot{\tilde{R}}_{i}(t,k)(\tilde{R}_i(t,k))^{-1} \\
\nonumber
& = \sum_{j \in \mathcal{N}_i}a_{ij}(\tilde{Q}_{j}(t,k)\tilde R_{ji}(t,k) - \tilde{Q}_i(t,k)),
\end{align}
or
\begin{align}
\label{eq:niss:2}
& \quad ~\tilde{U}(t,k) + [(\dot{\tilde{R}}_i(t,k)(\tilde{R}_i(t,k))^{-1})^T, 0 ]^T \\
\nonumber
& = \sum_{j \in \mathcal{N}_i}a_{ij}(\tilde{Q}_{ij}(t,k)\tilde R_{ji}(t,k) - [I_k, 0]^T),
\end{align}
where $\tilde{Q}_{ij}(t,k) = \tilde{Q}_i^T(t,d)\tilde{Q}_j(t,k)$ and $\tilde R_{ji}(t,k) = \tilde R_j(t,k) (\tilde R_i(t,k))^{-1}$. 
In the following we take equation \eqref{eq:niss:2} as the starting point and define 
\begin{align*}
\tilde{V}_i(t,k) & = \sum_{j \in \mathcal{N}_i}a_{ij}(\tilde{Q}_{ij}(t,k)\tilde{R}_{ji}(t,k) - [I_k, 0]^T),
\end{align*}
which is equivalent to \eqref{eq:main:1}.
Now, since $\dot{\tilde{R}}_i(t,k)(\tilde{R}_i(t,k))^{-1}$ is upper triangular we can obtain the lower part of $\tilde{U}(t,k)$ as 
$\text{low}(\tilde{V}_i)$. Since $\tilde{U}(t,d)$ is an element of $\mathsf{so}(d)$ for all $t$, it holds that the upper part of $\tilde{U}_i(t,k)$ is given by $[I_k, 0]^T(\text{low}(V_i(t,k))^T[I_k, 0]^T)$. Thus,
\begin{align*}
\tilde U_i(t,d) & = [\text{low}(\tilde V_i(t,k)), 0] - [\text{low}(\tilde V_i(t,k)), 0]^T, \\
\tilde U_i(t,k) & = \text{low}(\tilde V_i(t,k)) - [I_k, 0]^T(\text{low}(\tilde V_i(t,k))^T[I_k, 0]^T),
\end{align*}
which is equivalent to \eqref{eq:main:4} and \eqref{eq:main:2}, respectively. Now we can solve for $\dot{\tilde{R}}_i(t,k)$ in $\tilde{V}_i(t,k)$, $\tilde{U}_i(t,k)$, and $\tilde{R}_i(t,k)$. The solution is 
$$\dot{\tilde R}_i(t,k) = \text{up}((\tilde V_i(t,k) - \tilde U_i(t,k))\tilde R_i(t,k)),$$
which is equivalent to \eqref{eq:main:3}.

Since $h$ is a diffeomorphism, it holds that when the $Z_i(t,k)$'s converge to a full rank matrix $\bar{Z}$ as $t$ goes to infinity, the $\tilde Q_i(t,k)$'s converge to $\bar{Q}$ and the $\tilde{R}_i(t,k)$'s converge to $\bar{R}$ as $t$ goes to infinity, where $[\bar{Q}, \bar{R}]$ is the unique QR-decomposition of $\bar{Z}$ with positive diagonal elements in $\bar R$. It indeed holds that the $Z_i(t,k)$'s converge to a full rank matrix $\bar{Z}$ as $t$ goes to infinity, since the assumption is that $Z(0,t)$ is in $\mathcal{D}_{Z, \text{full}}(k)$.

\textbf{2)} Here we show that $h^{-1}(\mathcal{D}_{QR, \text{full}}(k)) \subset \mathcal{D}_{Z, \text{full}}(k)$.
Let $[Q^T(0,k) ,R^T(0,k)]^T$ in $\mathcal{D}_{QR, \text{full}}(k)$ be the initial condition for the dynamical system governed by \eqref{eq:main:1}-\eqref{eq:main:2}. At time $t$, let $Z_i(t,k) = h([Q_i^T(t,k), R_i^T(t,k)]^T) = Q_i(t,k)R_i(t,k)$, where $Q_i(t,k)$ and $R_i(t,k)$ are the matrices corresponding to agent $i$ at time $t$.
Since $[Q^T(0,k) ,R^T(0,k)]^T$ is in $\mathcal{D}_{QR, \text{full}}(k)$ and $h$ is a diffeomorphism, it holds that $Z_i(t,k)$ has full rank for all $i$, $t$ and all the $Z_i(t,k)$'s converge to a full rank matrix $\bar{Z}$. What remains to be shown, is that $Z(t,k)$ has the dynamics \eqref{eq:13}.
For all $i$, it holds that 
\begin{align*}
& \quad \: Q_i(t,d)^T\dot{Z}_i(t,k)R_i(t,k)^{-1} \\
& =  Q_i(t,d)^T\dot{ {Q}}_i(t,k)  +  Q_i(t,d)^T{Q}_i(t,k)\dot{ {R}}_i(t,k)R_i(t,k)^{-1} \\
& =  {U}_i(t,k) +  Q_i(t,d)^T{Q}_i(t,k)\dot{ {R}}_i(t,k)R_i(t,k)^{-1} \\
& =  {U}_i(t,k) +  [(\dot{ {R}}_i(t,k)R_i(t,k)^{-1})^T, 0]^T.
\end{align*}
Now we substitute the expressions for $U(t,d)$ and $\dot{R}_i(t,k)$ into the right-hand side of the last equation above. We obtain the following.
\begin{align*}
& \quad \:  U_i(t,k) + ~[\text{up}( V_i(t,k) -  U_i(t,k))^T, 0]^T \\
& = \text{low}(U_i(t,k))+ [\text{up}( V_i(t,k))^T, 0]^T \\
& = \text{low}(\text{low}( V_i(t,k)) - [I_k, 0]^T(\text{low}(  V_i(t,k))^T[I_k, 0]^T)) \\
& \quad +  [\text{up}( V_i(t,k))^T, 0]^T \\
& = V_i(t,k).
\end{align*}
Now, by definition it holds that 
\begin{align*}
{V}_i(t,k) & = \sum_{j \in \mathcal{N}_i}a_{ij}({Q}_{ij}(t,k){R}_{ji}(t,k) - [I_k, 0]^T).
\end{align*}
By multiplication of $Q_i(t,d)$ from the left and $R_i(t,k)$ from the right, respectively, we obtain  that
\begin{align*}
\dot{Z}_i(t,k) & = Q_i(t,d){V}_i(t,k)R_i(t,k) \\
& = \sum_{j \in \mathcal{N}_i}a_{ij}({Q}_{j}(t,k){R}_{j}(t,k) - Q_i(t,k)R_i(t,k)) \\
& = \sum_{j \in \mathcal{N}_i}a_{ij}(Z_j(t,k) - Z_i(t,k)). 
\end{align*} 
\hfill $\blacksquare$

\begin{rem}
In the proof we did not address the issue that $Q_i(t,d)$ needs to be in $\mathsf{SO}(d)$. If $k = d$, there is not a bijection $h$ between $\mathcal{D}_{Z}(k)$ and $\mathcal{D}_{QR}(k)$. However, since we only consider $k$ to be at most $d-1$, the $d-k$ last columns of $Q_i(t,d)$ can be chosen such that $Q_i(t,d)$ is in $\mathsf{SO}(d)$.
\end{rem}

\begin{rem}
For the QR-factorization in the proof, we are not using the standard convention where the ``$Q$-matrix'' is in $\mathbb{R}^{d \times d}$ and the ``$R$-matrix'' is in $\mathbb{R}^{d \times k}$; such a factorization is not unique (but can be chosen as a smooth function of $t$). Instead we let the ``$Q$-matrix'' be in $\mathbb{R}^{d \times k}$ and the ``$R$-matrix'' be in $\mathbb{R}^{k \times k}$. Those latter matrices are unique.
\end{rem}

\begin{rem}
In part \textbf{1)} of the proof, we can replace the $0$'s in the expression for $\tilde{U}_i(t,d)$ with something nonzero. Thus the expression for $\tilde{U}_i(t,d)$ is not unique, whereas the expression for $\tilde{U}_i(t,k)$ is. 
\end{rem}

\section{Illustrative example}
\begin{figure}[!th]
\centering
\includegraphics[scale=0.24]{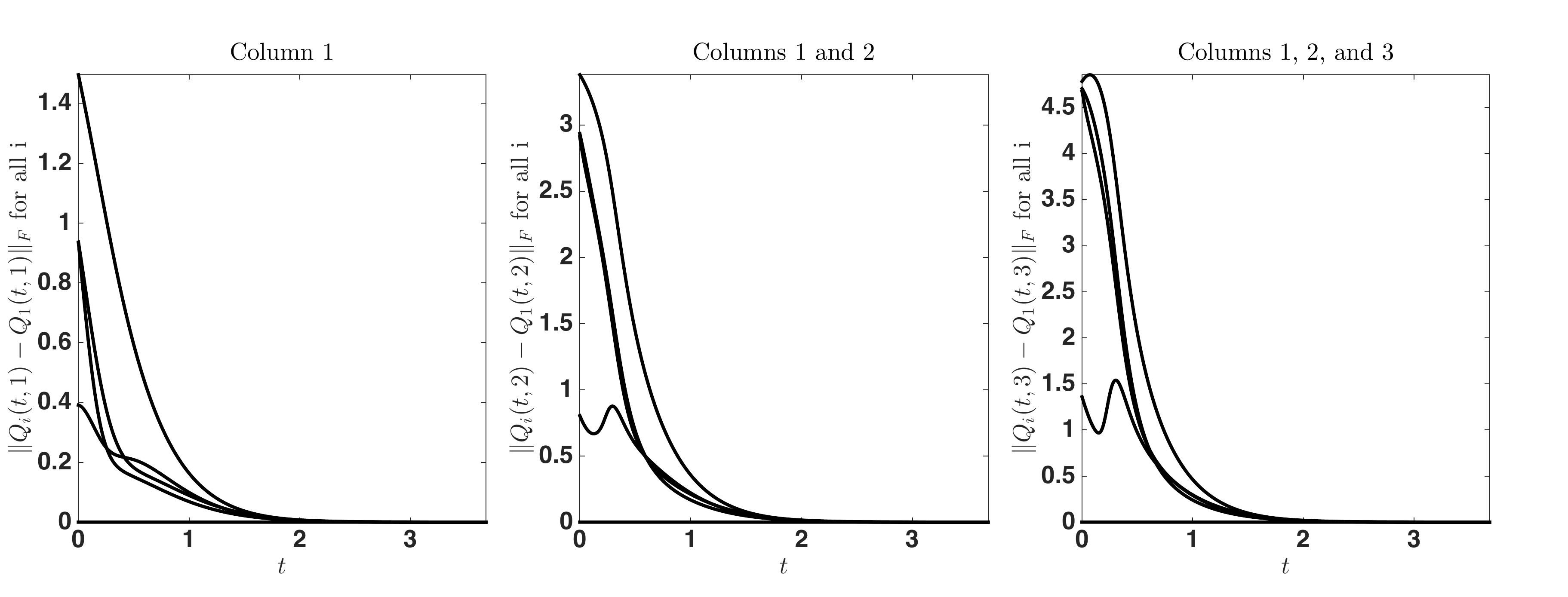}
\includegraphics[scale=0.24]{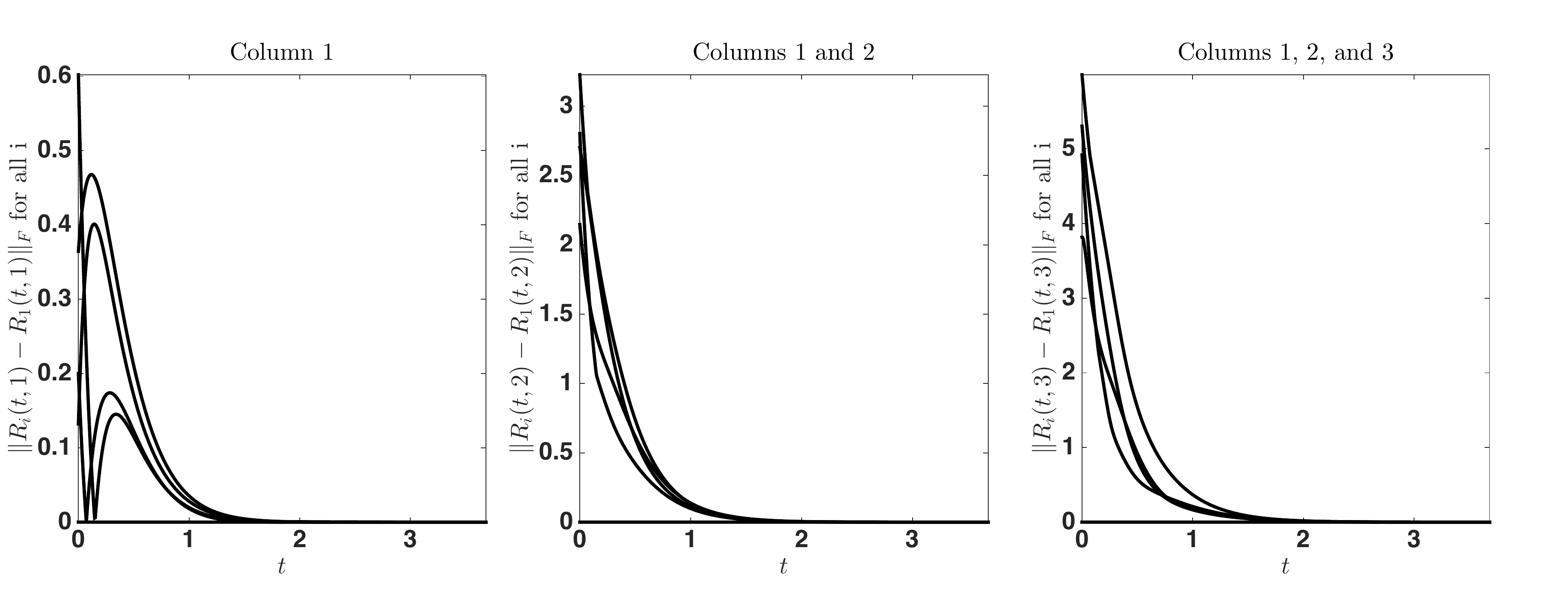}
\includegraphics[scale=0.24]{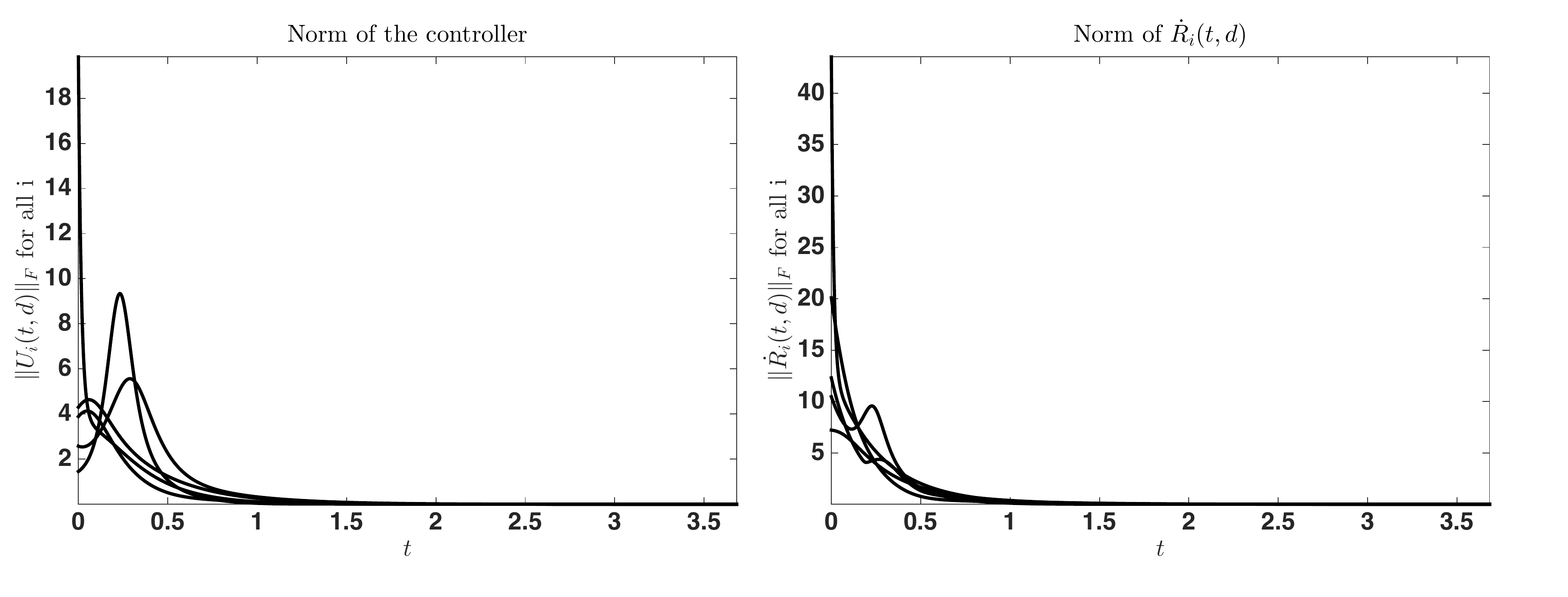}
\includegraphics[scale=0.23]{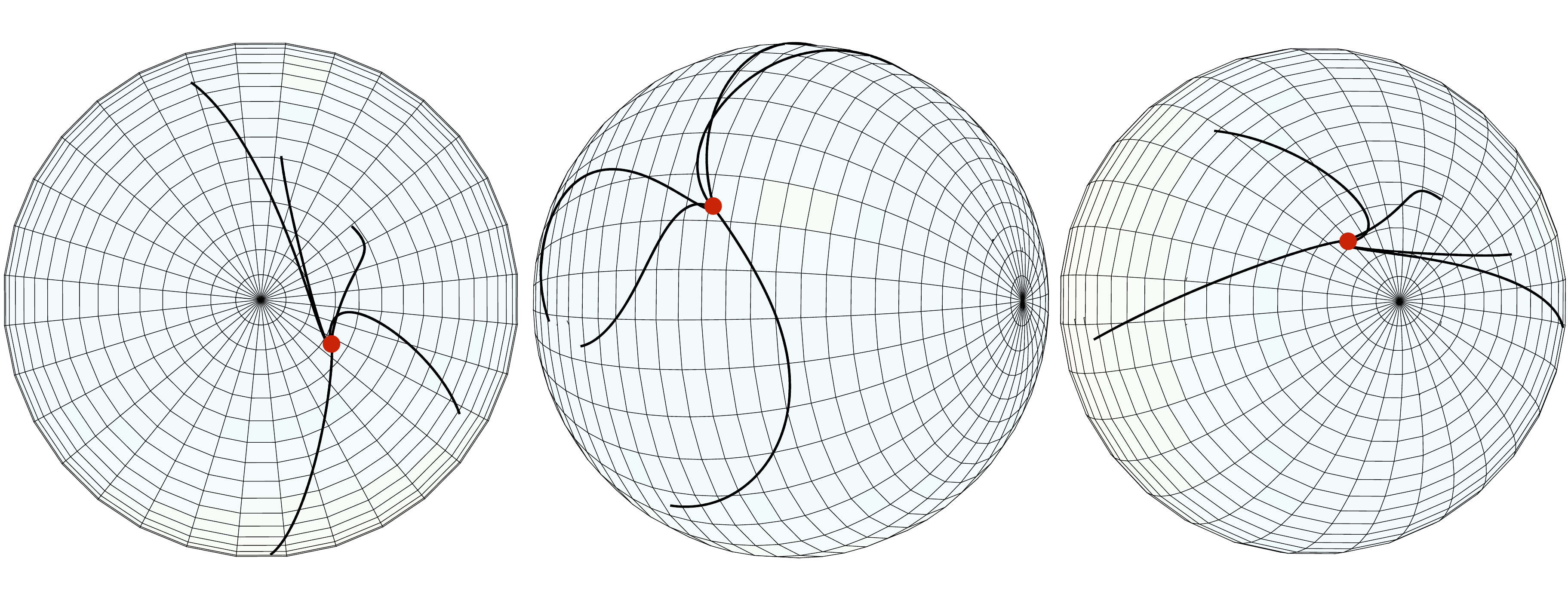}
\caption{Top three figures: convergence of $\|Q_i(t,k) - Q_1(t,k)\|_{\text{F}}$ to $0$ for all $i$ and $k$. Second from top three figures: convergence of $\|R_i(t,k) - R_1(t,k)\|_{\text{F}}$ to $0$ for all $i$ and $k$. Second from bottom two figures: convergence of $\|U_i(t,k)\|_{\text{F}}$ and $\|\dot{R}_i(t,k)\|_{\text{F}}$ to $0$ for all $i$. Bottom three figures: The convergence on the unit sphere of the individual columns of the $Q_i(t,d)$-matrices. Left: convergence of the first columns; middle: convergence of the second columns; right:  convergence of the third columns. The red disc denotes the unit vector that the columns converge to.}
 \label{fig:1}
\end{figure}

Fig.~\ref{fig:1} provides an illustrative example of the convergence when the controllers are given by Algorithm 1. 
Five agents were considered, i.e., $n = 5$, and the dimension $d$ was chosen to $3$. The topology of $\mathcal{G}$ was quasi-strongly connected but not strongly connected. The $Q_i(0,d)$'s and the $R_i(0,d)$'s were generated by first creating square matrices where each element was drawn from the Gaussian distribution with mean $0$ and variance $1$, and then performing QR-factorization of these matrices. The $a_{ij}$-scalars were generated by drawing samples from the uniform distribution with $(0,1)$ as support.
The simulations were conducted in Matlab by calling the function ode45.

\section{Conclusions}
We addressed the problem of column synchronization for rotation matrices in multi-agent systems. In this problem the first $k$ columns of the agents' rotation matrices shall be synchronized by only  using the corresponding $k$ first columns of the relative rotations between the agents.
For the proposed control scheme, we prove almost global convergence under quasi-strong interaction connectivity. The control scheme is based on the introduction of auxiliary state variables combined with a QR-factorization approach. The QR-factorization approach allows to separate the dynamics for the columns of interest from the remaining ones. This separation, in turn, enables us to only use the first $k$  columns of the relative rotations in the control design.

\bibliographystyle{agsm}       

\bibliography{ref.bib}

\end{document}